\documentclass[12pt,draft,noamsfonts]{amsart}

\usepackage{pb-diagram}

\newtheorem{thm}{Theorem}[section]   
\newtheorem{lem}[thm]{Lemma}         
\newtheorem{prop}[thm]{Proposition}  

\theoremstyle{definition}

\newtheorem{defn}[thm]{Definition}   

\theoremstyle{remark}

\newtheorem{rem}[thm]{Remark}        
\newtheorem{ex}[thm]{Example}        

\numberwithin{equation}{section}     

\newcommand{\secref}[1]{Section~\ref{#1}}
\newcommand{\thmref}[1]{Theorem~\ref{#1}}

\newcommand{\lemref}[1]{Lemma~\ref{#1}}
\newcommand{\propref}[1]{Proposition~\ref{#1}}
\newcommand{\defnref}[1]{Definition~\ref{#1}}

\newcommand{\exref}[1]{Example~\ref{#1}}

\newcommand{\Ind}{\ind}

\newcommand{\clsp}{\overline{\operatorname{\rm span}}}
\newcommand{\Ad}{\ad}

\newcommand{\Hom}{\operatorname{Hom}}

\newcommand{\rev}[1]{\widetilde{#1}}

\renewcommand{\c}[1]{\mathcal #1}

\renewcommand{\H}{\mathcal H}

\renewcommand{\L}{\mathcal L}

\newcommand{\cstar}{\ensuremath{C^*}-}

\newcommand{\spacetext}[1]{\quad\text{#1}\quad}

\newcommand{\midtext}[1]{\quad\text{#1}\quad}

\renewcommand{\phi}{\varphi}

\newcommand{\lip}[3]{
  {\vphantom\langle}_{#1}\!\!\left\langle{#2},{#3}\right\rangle  }
\newcommand{\rip}[3]{
  \left\langle{#2},{#3}\right\rangle_{\!{#1}}  }

\def\range{\operatorname{range}} 

\renewcommand{\H}{\mathcal H}
\newcommand{\A}{\mathcal A}
\newcommand{\B}{\mathcal B}
\newcommand{\C}{\mathcal C}
\newcommand{\K}{\mathcal K}

\DeclareMathOperator{\ad}{Ad}

\DeclareMathOperator{\id}{id}

\DeclareMathOperator{\ind}{Ind}
\DeclareMathOperator{\dashind}{-Ind}

\newcommand{\note}[1]{\footnote{#1}}

\begin{document}

\title{Naturality and Induced Representations}

\author{Siegfried Echterhoff}
\address{Fachbereich Mathematik und Informatik\\University of
M\"unster\\ 48149 M\"unster\\Germany}
\email{echters@uni-muenster.de}

\author{S. Kaliszewski}
\address{Department of Mathematics\\Arizona State University\\
 Tempe, Arizona 85287}
\email{kaz@math.la.asu.edu}

\author{John Quigg}
\address{Department of Mathematics\\Arizona State University\\
 Tempe, Arizona 85287}
\email{quigg@math.la.asu.edu}

\author{Iain Raeburn}
\address{Department of Mathematics\\University of Newcastle\\
 New South Wales 2308\\Australia}
\email{iain@maths.newcastle.edu.au}

\thanks{This research was partially supported by the 
Australian Research Council, and the Office of the Vice Provost for
Research at Arizona State University.}



\date{July 27, 1999; revised September 27, 1999}

\begin{abstract}
We show that induction of covariant representations 
for $C^*$-dynamical systems is natural in the sense that
it gives a natural transformation between certain crossed-product
functors.  This involves setting up
suitable categories of $C^*$-algebras and dynamical systems, and
extending the usual constructions of crossed products to 
define the appropriate functors.  {F}rom this point of view, 
Green's Imprimitivity
Theorem identifies the functors for which induction is a natural
equivalence.  Various spcecial cases of these results have previously
been obtained on an {\em ad hoc} basis. 

\end{abstract}

\maketitle


\section{Introduction}
\label{intro}

Induced representations and the imprimitivity theorems which
characterize them are a fundamental tool in the representation theory of
dynamical systems and crossed products. In the powerful formulation
of Rieffel and Green, induction is done by tensoring with a
Hilbert bimodule, and the imprimitivity theorem tells us how to expand
the left action to make this bimodule an imprimitivity bimodule.  In
recent applications, it has been necessary to know that this induction
process is compatible with other constructions involving crossed
products.
Verifying such
compatibility can be painful: it is often obvious that
everything must work because ``induction is natural'', yet
technically hard to sort out the details. And afterwards one is left
feeling that one must have missed the point: the techniques are vaguely
familiar even if the particular application isn't.

We have found that it is more
satisfactory to phrase our questions and results 
directly in terms of Hilbert bimodules,
rather than in terms of induced representations
themselves. Thus in \cite{KQR-DR} we viewed a Hilbert $A$--$B$ bimodule
$X$ as an arrow from $A$ to $B$, and said that a diagram 
\begin{equation}\label{testq}
\begin{diagram}
\node{A}
        \arrow{s,l}{Y}
        \arrow[2]{e,t}{X}
\node[2]{B}
        \arrow{s,r}{W}\\
\node{C}
        \arrow[2]{e,t}{Z}
\node[2]{D}
\end{diagram}
\end{equation} 
of Hilbert bimodules ``commutes in the strong sense'' if $Y\otimes_C
Z\cong X\otimes_B W$ as $A$--$D$ bimodules; if so, the
Rieffel induction  processes
$\pi\mapsto X\text{-}\Ind(W\text{-}\Ind\pi)$ and 
$\pi\mapsto Y\text{-}\Ind(Z\text{-}\Ind\pi)$
yield representations which
are systematically equivalent in a way compatible with intertwining
maps, direct-sum decompositions and continuity.

Here we make more precise the idea that we have been verifying the
naturality of induction. Category theory tells us exactly what
mathematicians should mean when they talk about naturality, and how we
need to set things up to make sense of it. Here we  build a category in
which the objects are
$C^*$-algebras, in which the morphisms are
given by Hilbert bimodules, and in which the commutativity of a diagram like
(\ref{testq}) says precisely that 
$Y\otimes_C Z\cong X\otimes_B W$. 
We then apply this framework to the induction process arising from
Green's imprimitivity theorem \cite[Proposition~3]{GreLS}, proving that
that theorem actually gives a natural transformation between 
certain crossed-product
functors which take values in this category.
This theorem includes as special cases a number of results in the
recent literature which could be paraphrased as saying ``Green
induction is compatible with Morita equivalence''.

Part of our point here is that taking the trouble to formulate theorems
in this categorical framework will result in a theory which is more
robust and more directly applicable.  
However, many applications require us to consider coactions or twisted
crossed products, and setting up the appropriate categories, functors,
and naturality theorems for these situations 
can take a good deal of technical effort. 
So in order to illustrate our approach, we 
concentrate here on one naturality theorem
which requires only standard techniques, and leave equivariant
versions involving coactions and twists to a more comprehensive sequel.  
We have tried hard to make this paper accessible to anyone
familiar with the basic material on Hilbert modules and Morita
equivalence implicit in the formulation of Green's theorem; 
this can be found, for example, in the early chapters of
\cite{RW}.  

\medskip

We begin in \secref{C}
by discussing our category $\mathcal{C}$.
The objects in $\mathcal{C}$ will be $C^*$-algebras, 
a morphism from one
$C^*$-algebra
$A$ to another $B$ will be given by 
a Hilbert $A$--$B$ bimodule $X$, 
and the composition of morphisms will be given by the balanced
tensor product of bimodules. There is asymmetry here: $X$ will be a right
Hilbert $B$-module, and $A$ will act by adjointable operators on $X$, in
the sense that there is a nondegenerate homomorphism $\kappa:A\to
\L(X_B)$ describing the action. We call such a bimodule $X$ a
right-Hilbert bimodule to stress that the Hilbert module structure is on
the right.\note{These and
related objects have been called {\em Hilbert bimodules} 
(for example, in \cite{PimCC}), or
{\em $C^*$-correspondences} (for example, in \cite{MS-TA}), but we
have preferred to stick with the name used in
\cite{KQR-DR} to stress that the Hilbert-module structure is on the
right.} The morphisms will actually be isomorphism classes $[_AX_B]$
of these right-Hilbert bimodules; this is necessary, because, for
example,  we expect $_BB_B$ to give the identity morphism on $B$, so
${}_AX_B$ and 
${}_A(X\otimes_B B)_B$ should define the same morphism in $\C$.
In \secref{C} we also show that the
isomorphisms in $\mathcal{C}$ are exactly the (isomorphism classes of)
imprimitivity bimodules, and prove that every morphism $[_AX_B]$ 
in $\C$ can be factored as $[{}_CY_B]\circ[\phi]$, 
where ${}_CY_B$ is an imprimitivity bimodule and 
$\phi:A\to M(C)$ is a nondegenerate homomorphism. 

In \secref{funct-sec} we 
show how to view the crossed-product construction as a
functor. We first have to add actions to the objects and morphisms of
$\mathcal{C}$ to build a category $\mathcal{A}(G)$ whose objects are
dynamical systems associated with a fixed locally compact group $G$;
forming the crossed product then gives a functor
$(A,\alpha)\mapsto A\times_\alpha G$ from $\mathcal{A}(G)$ to
$\mathcal{C}$. Because the objects in $\mathcal{C}$ are $C^*$-algebras
rather than isomorphism classes of $C^*$-algebras, we have to be
careful to nominate a particular $C^*$-algebra as the crossed product.
We discuss our nominee and its relationship to the universally
defined crossed product which we usually prefer. 
Functoriality requires that we give a parallel construction of crossed
products of Hilbert modules; the construction
we have chosen is different from
that of \cite{ComCP} and \cite{CMW-CP} 
and may be of some independent interest.

Our main theorem is proved in \secref{nat}. 
Green's construction  associates to
each object $(A,\alpha)$ in $\mathcal{A}(G)$ and each closed subgroup
$H$ of $G$ a
right-Hilbert $(A\times_\alpha G)$--$(A\times_\alpha H)$ bimodule
$X_H^G(A,\alpha)$;
his
Imprimitivity
Theorem says
that
$X_H^G(A,\alpha)$
is also a $(C_0(G,A)\otimes_{\alpha\otimes\tau}G)$--$(A\times_\alpha H)$ 
imprimitivity bimodule.  Our theorem says that the
assignment
$(A,\alpha)\mapsto [X_H^G(A,\alpha)]$ is a natural transformation.  
To prove this, 
we have to show that every morphism in $\mathcal{A}(G)$
gives rise to a certain commutative diagram in $\mathcal{C}$. We do this by
factoring the morphism into an ordinary homomorphism and an
imprimitivity bimodule; handling the homomorphism is straightforward,
and we deal with the imprimitivity bimodule by adapting a powerful
linking-algebra technique from \cite[\S4]{ER-ST}. We find it intriguing
that ideas developed to meet the demands of nonabelian duality are now
feeding back into the theory of ordinary crossed products --- both in
the technical sense, as in our use of the linking-algebra technique, and
in the motivational sense, in that we were led to formulate our theorem
through our attempts to handle more complicated problems
involving coactions.
 

\section{The Category $\c C$}
\label{C}

Let $B$ be a $C^*$-algebra.  A {\em Hilbert $B$-module\/}
is a vector space
$X$ which is a right $B$-module equipped with a positive definite
$B$-valued sesquilinear form $\rip{B}{\cdot}{\cdot}$ satisfying
\begin{equation}\label{B-mod-eq}
\rip{B}{x}{y\cdot b} = \rip{B}{x}{y}b\spacetext{and}
\rip{B}{x}{y}^* = \rip{B}{y}{x}
\end{equation}
for all $x,y\in X$ and $b\in B$, 
and which is complete in the norm $\|x\| =
\|\rip{B}{x}{x}\|^{\frac12}$.  
For the general theory of Hilbert modules, we refer to 
\cite{LanHC} or \cite{RW}.  

\begin{defn}\label{rHb-defn}
Let $A$ and $B$ be $C^*$-algebras, and let $X_B$ be a Hilbert
module which is {\em full} in the sense that
$\overline{\rip{B}{X}{X}}=B$.  We say that $X$ is  
a {\em right-Hilbert $A$--$B$ bimodule\/} if it is  
a nondegenerate left $A$-module satisfying
\begin{equation}\label{rHb-eq}
a\cdot(x\cdot b) = (a\cdot x)\cdot b\spacetext{and} 
\rip{B}{a\cdot x}{y} = \rip{B}{x}{a^*\cdot y}
\end{equation}
for all $a\in A$, $x,y\in X$, and $b\in B$.
(This is equivalent to having a nondegenerate homomorphism $\kappa\colon
A\to \L(X_B)$ and putting $a\cdot x = \kappa(a)x$.)   
An {\em isomorphism} of right-Hilbert bimodules 
is a bijective linear map $\Phi\colon {}_AX_B
\to {}_AY_B$ such that 
\begin{itemize}
\item[(i)] $\Phi(a\cdot x) = a\cdot\Phi(x)$,
\item[(ii)] $\Phi(x\cdot b) = \Phi(x)\cdot b$, and
\item[(iii)] $\rip{B}{\Phi(x)}{\Phi(y)} = \rip{B}{x}{y}$.
\end{itemize}
(Property (ii) is redundant ---  it follows from 
(iii) by considering $\|\Phi(x\cdot b)-\Phi(x)\cdot b\|^2$
--- but we include it to emphasize that $\Phi$ is a bimodule
homomorphism.)
\end{defn}

\begin{ex}\label{rHb-ex} 
If $B$ is a $C^*$-algebra, then $B$
becomes a 
right-Hilbert $B$--$B$ bimodule by putting
$$a\cdot b\cdot c = abc\spacetext{and} \rip{B}{a}{b} = a^*b$$ 
for $a,b,c\in B$.  
More generally, if $\phi\colon A\to M(B)$ is a
nondegenerate homomorphism, then $B_B$ becomes a
right-Hilbert $A$--$B$ bimodule with left action given by 
$$a\cdot b = \phi(a)b.$$
\end{ex}

The morphisms from $A$ to $B$ will be the isomorphism classes of
right-Hilbert
$A$--$B$ bimodules (denoted with square brackets);  we need to pass to
isomorphism classes to show that  composition of morphisms has the
necessary properties.    Two ordinary homomorphisms may give the same
morphism in $\C$, but only if they differ by an inner automorphism:

\begin{prop}
Let $\phi$ and $\psi$ be nondegenerate homomorphisms of $A$ into
$M(C)$.
Then $[\phi]=[\psi]$ in $\C$ if and only if there exists $u\in
UM(C)$
such that $\psi=\Ad u\circ\phi$.
\end{prop}

\begin{proof}
If $\psi=\Ad u\circ\phi$, then the
map $c\mapsto uc$ is Hilbert module 
automorphism of $C_C$ which intertwines
the left actions coming from $\phi$ and $\psi$.

For the converse, suppose $[\phi]=[\psi]$ in $\C$, so there
exists a
linear bijection $ L\colon C\to C$ such that
\[
 L(cd) =  L(c)d,\qquad
 L(c)^* L(d) = c^*d,\midtext{and}
 L(\phi(a)c) = \psi(a) L(c)
\]
for each $c,d\in C$.  Define $R\colon C\to C$ by
$R(c) =  L^{-1}(c^*)^*$.  Then the first two of these properties
imply that $( L,R)$ is an invertible double centralizer of $C$,
so there exists an invertible element $u$ of $M(C)$ such that
$ L(c) = uc$ for all $c$.  Since 
\[
u^{-1}c =  L^{-1}(c) =  L^{-1}(c^{**})^{**} = 
R(c^*)^* = (c^*u)^* = u^*c,
\]
it follows that $u$ is unitary, and $ L(\phi(a)c) =
\psi(a) L(c)$ implies that $\psi=\Ad u\circ \phi$.  
\end{proof}

To define composition we use the internal tensor product of Hilbert
modules.  
Let ${}_AX_B$ and ${}_BY_C$ be right-Hilbert bimodules.  Then the
algebraic tensor product $X\odot Y$ becomes a not-necessarily-complete
right-Hilbert $A$--$C$ bimodule with actions given by
$$a\cdot (x\otimes y) = a\cdot x\otimes y\spacetext{and} (x\otimes
y)\cdot b = x\otimes y\cdot b$$
and a not-necessarily-definite $C$-valued inner product given by
$$\rip{C}{x\otimes y}{z\otimes w} = \rip{C}{y}{\rip{B}{x}{z}\cdot w}$$
\cite[Proposition~4.5]{LanHC}.
Modding out by the elements of length $0$ and completing
gives a right-Hilbert $A$--$C$ bimodule $X\otimes_B Y$
in which $x\cdot b\otimes y = x\otimes b\cdot y$.  We emphasize
that it is not obvious that the action of $a\in A$ on $X\odot Y$
extends to $X\otimes_B Y$; this requires a complete positivity argument
as in \cite[page~42]{LanHC}.
Note that $X\otimes_BY$ is full, since $X$ and $Y$ are:
$$\overline{\rip{C}{X\otimes_BY}{X\otimes_BY}}
= \overline{\rip{C}{Y}{\rip{B}{X}{X}\cdot Y}}
= \overline{\rip{C}{Y}{B\cdot Y}}
= \overline{\rip{C}{Y}{Y}} = C.$$
We call 
$X\otimes_B Y$ the {\em balanced\/} or {\em internal\/}
{\em tensor product\/}.

\begin{prop}\label{C-cat}
There is a category $\C$ in which the objects are \cstar algebras, and
in which the morphisms from $A$ to $B$  are the
isomorphism classes of (full) right-Hilbert $A$--$B$ bimodules.
The composition of $[X]\colon A\to B$ with $[Y]\colon B\to C$
is the isomorphism class of the internal tensor product bimodule
${}_A(X\otimes_B Y)_C$;  the identity morphism 
on $A$ is the isomorphism class
of the right-Hilbert bimodule ${}_AA_A$.
\end{prop}

\begin{rem}
In any category, $\Hom(A,B)$ has to be a set for each pair of
objects $A$ and $B$.  This is not obviously true in $\C$ 
unless we limit the size of the bimodules involved. 
We can do this by 
considering only $C^*$-algebras and Hilbert
modules with dense subsets whose cardinalities do not exceed a fixed
large cardinal. 
For example, we could consider only separable $C^*$-algebras and
bimodules.  
In practice, these issues should never present an unassailable problem,
and  we shall ignore them.
\end{rem}

\begin{proof}[Proof of Proposition~\ref{C-cat}]
We first claim that the composition of morphisms is well-defined.
Suppose we
have right-Hilbert bimodule isomorphisms
$\Phi\colon{}_AX_B\to{}_AZ_B$ and
$\Psi\colon{}_BY_C\to{}_BW_C$.
Then $\Phi\otimes\Psi\colon X\odot Y \to Z\otimes_BW$ is easily seen to
preserve the actions and inner product, so extends to an isometric
bimodule map $\Phi\otimes_B\Psi$.  The map
$\Phi^{-1}\otimes_B\Psi^{-1}$ is an inverse for  $\Phi\otimes_B\Psi$,
which is therefore a right-Hilbert $A$--$C$
bimodule isomorphism of $X\otimes_BY$ onto $Z\otimes_BW$.

Next we establish that composition of morphisms in $\C$ is associative;
  it suffices to show that $X\otimes_B(Y\otimes_CZ)$ and
$(X\otimes_BY)\otimes_CZ$ are isomorphic for any right-Hilbert
bimodules
${}_AX_B$, ${}_BY_C$, and ${}_CZ_D$.  But straightforward calculations
show
that the usual linear isomorphism of $X\odot(Y\odot Z)$ onto $(X\odot
Y)\odot
Z$ respects the module actions and right inner products, so extends to
the
desired isomorphism.

Finally, note that the maps $a\otimes x \mapsto a\cdot x$ and $y\otimes
a\mapsto y\cdot a$ extend to isomorphisms $A\otimes_AX\cong X$ and
$Y\otimes_AA\cong Y$ for any right-Hilbert bimodules ${}_AX_B$ and
${}_BY_A$.
Hence 
$[{}_AA_A]$
is an identity morphism from $A$ to $A$ in $\c C$.  
\end{proof}

It is implicit in Connes' definition  of Morita equivalence
(\cite[p.~155]{ConNG})
that a morphism $[X]$ is invertible in $\C$ precisely when $X$ is
an imprimitivity bimodule.  
This seemed obvious to us at first 
but we found it surprisingly hard to prove. 

By saying that a right-Hilbert bimodule ${}_AX_B$ is an imprimitivity
bimodule we mean that $X$ is also a (full) left-Hilbert $A$-module
in such a way that 
\[
\lip{A}{x}{y}\cdot z = x\cdot\rip{B}{y}{z}\midtext{and}
\lip{A}{x\cdot b}{y} = \lip{A}{x}{y\cdot b^*}.
\]
Equivalently, ${}_AX_B$ is an imprimitivity bimodule if and only if
the canonical map $\kappa\colon A\to\L(X_B)$ induced by the left action
of $A$ on $X_B$ is an isomorphism of $A$ onto the algebra $\K(X_B)$ of
compact operators on $X_B$ \cite[Proposition~3.8]{RW}.   

\begin{prop}\label{ib-isom-prop}
Suppose ${}_AX_B$ is a right-Hilbert bimodule.  Then $[X]$ 
is an isomorphism in the category $\C$ if and only if $X$ is an 
$A$--$B$ imprimitivity bimodule. 
\end{prop}

One direction of \propref{ib-isom-prop} is easy: if ${}_AX_B$ is an
imprimitivity bimodule, and ${}_B\rev{X}_A$ is its reverse bimodule,
then the maps $x\otimes\rev{y}\mapsto\lip{A}{x}{y}$ and $\rev{x}\otimes
y\mapsto\rip{B}{x}{y}$ are isomorphisms of $X\otimes_B\rev{X}$ onto
${}_AA_A$ and $\rev{X}\otimes_AX$ onto ${}_BB_B$, respectively (see
\cite[Proposition~3.28]{RW}), and hence $[\rev{X}]$ is an inverse for
$[X]$ in $\C$.  

For the other direction we want to use a representation-theoretic
argument, so we briefly discuss how right-Hilbert bimodules can be
represented by bounded operators.  These ideas are fairly standard: the
analogues for imprimitivity bimodules appear in \cite[Section~2]{ER-MI}.  

\begin{defn}
\label{rep}
A {\em representation} $(\pi_A,\pi_X,\pi_B)$
of a right-Hilbert bimodule ${}_AX_B$ is a
triple consisting of nondegenerate
representations $\pi_A$ and $\pi_B$ of $A$ and $B$ on Hilbert spaces
$\H_A$ and $\H_B$, respectively, together with a linear map
$\pi_X\colon X\to \B(\H_B,\H_A)$ such that
\begin{itemize}
\item[(i)]
$\pi_X(a\cdot x) = \pi_A(a)\pi_X(x)$,
\item[(ii)]
$\pi_X(x\cdot b) = \pi_X(x)\pi_B(b)$,\midtext{and}
\item[(iii)]
$\pi_B(\rip{B}{x}{y}) = \pi_X(x)^*\pi_X(y)$
\end{itemize}
for all $a\in A$, $x,y\in X$, and $b\in B$.  
(Again, property (ii) follows from property (iii), but we include it
so the definition explicitly says that $\pi$ preserves the Hilbert-module
structure of $X$.)  
\end{defn}

\begin{ex}\label{rep-ex}
Suppose ${}_AX_B$ is a right-Hilbert bimodule and $\pi_B$ is any
nondegenerate representation of $B$ on a Hilbert space $\H_B$.  The
{\em induced representation} $\pi_A=X\dashind\pi_B$ of $A$ acts in
$\H_A=X\otimes_B\H_B$ according to $\pi_A(a)(x\otimes\xi) = (a\cdot
x)\otimes\xi$ \cite[Proposition~2.66]{RW}.  
We claim that there are bounded operators
$\pi_X(x)\in\L(\H_B,\H_A)$ such that $\pi_X(x)\xi = x\otimes\xi$, and
that $(\pi_A,\pi_X,\pi_B)$ is then a nondegenerate representation of
${}_AX_B$ on $(\H_A,\H_B)$.  Indeed, it is easy to check that
$\pi_X(x)$ is bounded:
\[
\|\pi_X(x)\xi\|^2 = (x\otimes\xi\mid x\otimes\xi) =
(\pi_B(\rip{B}{x}{x})\xi\mid\xi) \leq
\|\pi_B(\rip{B}{x}{x})\|\,\|\xi\|^2;
\]
the relation $\pi_A(a)\pi_X(x) = \pi_X(a\cdot x)$ is immediate from the
definition of $\pi_A$; the relation $\pi_X(x)\pi_B(b) = \pi_X(x\cdot
b)$ follows from the balancing property $(x\cdot b)\otimes\xi =
x\otimes\pi_B(b)\xi$ of the tensor product; and the definition of the inner
product on $X\otimes_B\H_B$ implies that $\pi_B(\rip{B}{x}{y}) =
\pi_X(x)^*\pi_X(y)$.  
\end{ex}

In this example, the representation $\pi_A$ is by
definition the restriction of a representation $\pi_\L$ of the algebra
$\L(X_B)$ of adjointable operators on $X$ (strictly speaking, $\pi_A =
\pi_\L\circ\kappa$).  The following observation will be useful in the
proof of \propref{ib-isom-prop}. 

\begin{lem}
Let $X_B$ be a right Hilbert module, let $\pi_B$ be a faithful
representation of $B$ on $\H_B$, and consider the representation
$(\pi_\L, \pi_X, \pi_B)$ of the right-Hilbert bimodule ${}_{\L(X)}X_B$
discussed in \exref{rep-ex}.  Then $\pi_\L$ is faithful, and
\[
\pi_\L(\K(X_B)) = \pi_X(X)\pi_X(X)^* := \clsp\{\pi_X(x)\pi_X(y)^*\mid
x,y\in X\}.
\]
\end{lem}

\begin{proof}
By definition, $\K(X_B)$ is the closed span of the rank-one operators
$\{\Theta_{x,y}\mid x,y\in X\}$ \cite[Definition~2.24]{RW}, so it is
enough to show that $\pi_\L(\Theta_{x,y}) = \pi_X(x)\pi_X(y)^*$.  We
can verify this by applying both sides to vectors of the form
$x\otimes\xi$, which densely span $X\otimes_B\H_B$.  
We then have
\begin{align*}
\pi_\L(\Theta_{x,y})(z\otimes\xi)
& =\Theta_{x,y}(z)\otimes\xi
=x\cdot\rip{B}{y}{z}\otimes\xi
= x\otimes\pi(\rip{B}{y}{z})\xi\\
& =x\otimes\pi_X(y)^*\pi_X(z)\xi
=\pi_X(x)(\pi_X(y)^*\pi_X(z)\xi)\\
&=(\pi_X(x)\pi_X(y)^*)(z\otimes\xi),
\end{align*}
as required.  

Since $\pi_B$ is faithful and $X$ is a $\K(X_B)$--$B$ imprimitivity
bimodule, the induced representation $\pi_\K$ is faithful (this follows
from the Rieffel correspondence \cite[Proposition~3.24]{RW}), and hence
so is its extension $\pi_\L$ to $\L(X_B) = M(\K(X_B))$.  
\end{proof}

\begin{proof}[Proof of \propref{ib-isom-prop}.]  
We have to show that if there is a right-Hilbert bimodule ${}_BY_A$
such that $X\otimes_BY\cong {}_AA_A$ and $Y\otimes_AX\cong {}_BB_B$ as
right-Hilbert bimodules, then $X$ is an $A$--$B$ imprimitivity
bimodule (and in fact $Y$ will then be isomorphic to $\rev{X}$).  We
choose a faithful nondegenerate representation $\pi_B$ of $B$ on $\H_B$,
and use the inducing construction of \exref{rep-ex} to give a
representation $(\pi_A, \pi_X, \pi_B)$ of ${}_AX_B$ on $(\H_B, \H_A)$
with $\H_A= X\otimes_B\H_B$.  We shall prove the result by showing that the
presence of $Y$ implies that $\pi_A$ is an isomorphism of $A$ onto the
image $\clsp\{\pi_X(x)\pi_X(y)^*\}$ of the algebra $\K(X_B)$ of compact
operators.  

We begin by noting that because $A$ acts faithfully on itself and
${}_AA_A\cong X\otimes_BY$, $A$ must act faithfully on $X$: if $a\cdot
x = 0$ for all $x$, then $a\cdot(x\otimes y) = 0$ for all $y$, whence
$ab=0$ for all $b\in A$, and $a=0$.  Thus $\kappa\colon A\to\L(X_B)$
is
injective, and the injectivity of $\pi_\L$ implies 
that $\pi_A=\pi_\L\circ\kappa$
is faithful.  To show that $\pi_A$ has the right image, we construct a
representation $(\pi_B, \sigma, \pi_A)$ of $Y$ on $(\H_A,\H_B)$.  

Let $\rho\colon Y\to
\B(\H_A,Y\otimes_A\H_A)$ be the representation induced from $\pi_A$,
given as usual by $\rho(y)\zeta =
y\otimes\zeta$.  By assumption, there exists a right-Hilbert bimodule
isomorphism $\Psi\colon Y\otimes_AX\to B$, which we use to define a
unitary $U\colon Y\otimes_A\H_A = Y\otimes_AX\otimes_B\H_B\to \H_B$ by
$$U(y\otimes x\otimes\xi) = \pi_B(\Psi(y\otimes x))\xi.$$
Now we define
$$\sigma(y) = U\circ\rho(y).$$
For $y\in Y$, $b\in B$, and $x\otimes\xi\in \H_A$ we get
\begin{align*}
\sigma(b\cdot y)(x\otimes\xi)
& = U(b\cdot y\otimes x\otimes\xi)= \pi_B(\Psi(b\cdot y\otimes x))\xi\\
& = \pi_B(b)\pi_B(\Psi(y\otimes x))\xi=\pi_B(b)U(y\otimes
x\otimes\xi)\\ 
& =\pi_B(b)\sigma(y)(x\otimes\xi),
\end{align*}
and we also have
$$\sigma(y)^*\sigma(z) = \rho(y)^*U^*U\rho(z)
= \rho(y)^*\rho(z) = \pi_A(\rip{A}{y}{z})$$
for all $y,z\in Y$, since $\rho$ has right coefficient map $\pi_A$.
Thus $(\pi_B,\sigma,\pi_A)$ is indeed a
right-Hilbert bimodule representation of $Y$.
In particular, we have $\pi_A(A) = \sigma(Y)^*\sigma(Y)$.

It only remains to show that $\sigma(Y)^* = \pi_X(X)$.
For this we
observe that for $y\in Y$, $x\in X$, and $\xi\in\H_B$ we have
$$\sigma(y)\pi_X(x)\xi = \sigma(y)(x\otimes\xi) = \pi_B(\Psi(y\otimes
x))\xi,$$
hence $\sigma(Y)\pi_X(X) = \pi_B(\Psi(Y\otimes_AX)) = \pi_B(B)$.  Using
$\pi_A(A) = \sigma(Y)^*\sigma(Y)$, it now follows that
\begin{gather*}
\pi_X(X) = \pi_X(A\cdot X) = \pi_A(A)\pi_X(X) =
\sigma(Y)^*\sigma(Y)\pi_X(X)\\
= \sigma(Y)^*\pi_B(B) = \sigma(Y\cdot B)^* = \sigma(Y)^*.
\end{gather*}
This completes the proof of \propref{ib-isom-prop}. 
\end{proof}

\begin{prop}\label{decomp}
For every morphism $[{}_AX_B]$ in $\C$, there exists a $C^*$-algebra
$C$, a nondegenerate homomorphism $\phi\colon A\to M(C)$, and an
imprimitivity bimodule ${}_CY_B$ such that
\[
[{}_AX_B] = [{}_CY_B]\circ[\phi]\midtext{in}\C.
\]
In other words, such that ${}_AX_B\cong {}_A(C\otimes_CY)_B$ as
a right-Hilbert bimodule. 
\end{prop}

\begin{proof}
Let $C=\K(X_B)$ and let $Y=X$ viewed as a $C$--$B$ imprimitivity
bimodule (see \cite[Proposition~3.8]{RW}). Since 
$\L(X_B)\cong M(\K(X_B))$ 
(\cite[Corollary~2.54]{RW}), we can view
the canonical nondegenerate
homomorphism $\kappa\colon A\to \L(X_B)$ as a nondegenerate
homomorphism $\phi:A\to M(\K(X_B))$.  
The map 
$c\otimes y\mapsto c\cdot y$ extends to an isomorphism 
of $C\otimes_C Y$ onto $X$.  
\end{proof}


\section{Crossed-Product Functors}
\label{funct-sec}

Our next goal is to
formalize the idea that assignments like $(A,\alpha)\mapsto
A\times_\alpha G$ are functors into $\C$; this requires that we
construct an equivariant category $\A(G)$ in which the objects are
dynamical systems $(A,G,\alpha)$.   

\begin{defn}\label{rHb-act-defn}
Let $G$ be a locally compact group, 
let ${}_AX_B$ be a right-Hilbert bimodule, and 
let $\alpha$ and $\beta$ be (strongly continuous) 
actions of $G$ on $A$ and
$B$.  An {\em $(\alpha,\beta)$-compatible action\/} of $G$ on $X$ is a 
homomorphism $\gamma$ of $G$ into the group  of 
invertible linear transformations on $X$
such that
\begin{itemize}
\item[(i)] $\gamma_s(a\cdot x)=\alpha_s(a)\cdot\gamma_s(x)$
\item[(ii)] $\gamma_s(x\cdot b)=\gamma_s(x)\cdot\beta_s(b)$
\item[(iii)] $\rip{B}{\gamma_s(x)}{\gamma_s(y)} =
\beta_s(\rip{B}{x}{y})$
\end{itemize}
for each $s\in G$, $a\in A$, $x,y\in X$, and $b\in B$;
and such that each map $s\mapsto\gamma_s(x)$ is continuous from $G$
into $X$.  
(As usual, (ii) follows from (iii).)  Two 
$(\alpha,\beta)$-compatible actions $\gamma$ and $\eta$ on $X$ and $Y$
are {\em equivariantly isomorphic} if there exists an isomorphism
$\Phi$ of $X$ onto $Y$ such that $\Phi(\gamma_s(x)) = \eta_s(\Phi(x))$
for each $s\in G$ and $x\in X$.  
\end{defn}

\begin{ex}
If $(A,\alpha)$ and $(B,\beta)$ are actions of $G$, and $\phi\colon
A\to M(B)$ is equivariant in the sense that $\phi(\alpha_s(a)) =
\beta_s(\phi(a))$, then $\beta$ is an $(\alpha,\beta)$
compatible action of $G$ on the bimodule ${}_AB_B$ of
\exref{rHb-ex}.  
\end{ex}

\begin{prop}\label{AG-cat}
Let $G$ be a locally compact group.  
\begin{itemize}
\item[(i)]
There is a category $\A(G)$ 
in which the objects are $C^*$-algebras with actions of $G$, and in
which the morphisms from $(A,\alpha)$ to $(B,\beta)$ are the
equivariant 
isomorphism classes of (full) right-Hilbert $A$--$B$ bimodules with
$(\alpha,\beta)$-compatible 
actions of $G$.  The composition of $[X,\gamma]\colon
(A,\alpha)\to (B,\beta)$ with $[Y,\eta]\colon (B,\beta)\to
(C,\epsilon)$
is the isomorphism class of the tensor product action
$(X\otimes_BY,\gamma\otimes_B\eta)$; the identity morphism on
$(A,\alpha)$ is $[{}_{(A,\alpha)}(A,\alpha)_{(A,\alpha)}]$.
\item[(ii)]
$[{}_{(A,\alpha)}(X,\gamma)_{(B,\beta)}]$
is an isomorphism in $\A(G)$ if and only if 
${}_AX_B$ is an imprimitivity bimodule and
$\alpha_s(\lip{A}{x}{y}) = \lip{A}{\gamma_s(x)}{\gamma_s(y)}$. 
\item[(iii)]
For every morphism $[{}_{(A,\alpha)}(X,\gamma)_{(B,\beta)}]$ 
in $\A(G)$, there exists an
isomorphism $[{}_{(C,\epsilon)}(Y,\eta)_{(B,\beta)}]$
and an $\alpha$--$\epsilon$ equivariant nondegenerate homomorphism
$\phi\colon A\to M(C)$ such that 
\[
[X,\gamma] = [Y,\eta]\circ[\phi,\epsilon]
\midtext{in} \A(G).
\]
\end{itemize}
\end{prop}

\begin{proof}
Adding actions to Propositions \ref{C-cat}, \ref{ib-isom-prop},
and~\ref{decomp} is routine, except possibly in statement (iii).  
In the proof of \propref{decomp}, we took $C=\K(X_B)$ and
$\phi=\kappa\colon A\to \L(X_B)=M(C)$, so we need to show that
$\gamma$ induces an action on $\K(X_B)$.  But for each $T\in
\L(X_B)$, $\epsilon_s(T)\colon x\mapsto
\gamma_s(T(\gamma_s^{-1}(x)))$ is an adjointable operator with
$\epsilon_s(T)^*=\epsilon_s(T^*)$, and in fact $T\mapsto
\epsilon_s(T)$ is an automorphism of $\L(X_B)$.  A quick
calculation shows that 
\begin{equation}
\label{cpts}
\epsilon_s(\lip{\K}{x}{y}) = \lip{\K}{\gamma_s(x)}{\gamma_s(y)},
\end{equation}
so $\epsilon_s$ restricts to an automorphism of $\K(X_B)$.  Since
$G$ acts continuously on $X$, $\epsilon$ gives an action of
$G$ on $\K$, and \eqref{cpts} says that
${}_{(C,\epsilon)}(X,\gamma)_{(B,\beta)}$ is an isomorphism in
$\A(G)$.  Now recall that $\kappa(\alpha_s(a))(x) =
\alpha_s(a)\cdot x = \gamma_s(a\cdot\gamma_s^{-1}(x))$ to see
that $\phi$ is $\alpha$--$\epsilon$ equivariant;  
it is then easy to check that the isomorphism $c\otimes y\mapsto
c\cdot y$ of $C\otimes_CY$ onto $X$ is
$\epsilon\otimes\eta$--$\gamma$ equivariant.  
\end{proof}

To construct the crossed product of a system $(A,G,\alpha)$, we begin
with the vector space $C_c(G,A)$ of continuous functions $f\colon G\to
A$ of compact support.  This is a $*$-algebra 
with the operations 
\[
f*g(s) = \int_G f(t)\alpha_t(g(t^{-1}s))\,dt\midtext{and}
f^*(s) = \alpha_s(f(s^{-1})^*)\Delta_G(s)^{-1}.
\]
A covariant representation of $(A,G,\alpha)$ on
$\H$ is a pair $(\pi,U)$ consisting of a nondegenerate representation
$\pi\colon A\to \B(\H)$ and a unitary representation $U\colon G\to
U(\H)$ such that $\pi(\alpha_s(a)) = U_s\pi(a)U_s^*$ for $s\in G$ and $a\in
A$.  Routine calculations show that there is a $*$-representation
$\pi\times U$ of $C_c(G,A)$ on $\H$, called the integrated form of
$(\pi,U)$, such that 
\[
\pi\times U(f)\xi = \int_G \pi(f(s))U_s\xi\,ds
\]
for all $\xi\in \H$.  (Inserting the vector $\xi$ is technically
helpful because it makes the integrand norm-continuous.)  
Since $\|\pi\times U(f)\|\leq \int_G \|f(s)\|\,ds$ for every $f$,
we can define a semi-norm on $C_c(G,A)$ by 
\note{
There is ostensibly a set-theoretic problem here since
it is not obvious that we are taking the supremum over a
well-defined set of real numbers.  But the numbers of the
form $\|\pi\times U(f)\|$ do form a subclass of $\bf R$, and every
subclass of a set is again a set (\cite[Proposition~4.6]{MenIM}).
}
\[
\|f\|_* = \sup\{\|\pi\times U(f)\|\colon (\pi,U)\text{\ is a
covariant representation of\ } (A,G,\alpha)\}; 
\]
it is a norm because the regular representation
$\tilde\pi\times \lambda$ on $L^2(G,\H_\pi)$ is faithful on
$C_c(G,A)$ if $\pi$ is faithful on $A$.  

We now define the crossed product $A\times_\alpha G$ to be the
completion of $C_c(G,A)$ in the norm $\|\cdot\|_*$.  (For this to be a
construction of a particular $C^*$-algebra rather than an isomorphism
class of $C^*$-algebras, we must be clear that forming completions is a
construction, but this can be achieved by defining 
the completion of a
normed space $X$ to be the closure of $X$ in its double dual.)
Notice that every representation $\pi\times U$ of $C_c(G,A)$ is
by definition continuous for $\|\cdot\|_*$, and hence extends to
a representation of $A\times_\alpha G$, which we continue to
denote by $\pi\times U$.    

\begin{rem}
In \cite{RaeOC}, a crossed product for the system $(A,G,\alpha)$ is any
$C^*$-algebra $B$ equipped with a nondegenerate homomorphism $i_A\colon
A\to M(B)$ and a strictly continuous homomorphism $i_G\colon G\to
UM(B)$ satisfying
\begin{itemize}
\item[(a)]
$(i_A,i_G)$ is covariant, in the sense that $i_A(\alpha_s(a)) =
i_G(s)i_A(a)i_G(s)^*$ for all $a\in A$ and $s\in G$;
\item[(b)]
if $(\pi,U)$ is a covariant representation, there is a nondegenerate
representation $\phi=\phi_{\pi,U}$ of $B$ such that 
$\pi=\phi\circ i_A$ and $U=\phi\circ i_G$;
and
\item[(c)]
the elements $i_A\times i_G(f) = \int_G i_A(f(s)) i_G(s)\,ds$ for $f\in
C_c(G,A)$ span a dense subspace of $B$. 
\end{itemize}
{F}rom our present perspective, we can view these axioms as
properties which allow us to identify isomorphic copies
of $A\times_\alpha G$.  To be more precise, if $(B,i_A,i_G)$ is a
crossed product in the sense of \cite{RaeOC}, then we claim that
$i_A\times i_G$ extends to an isomorphism of $A\times_\alpha G$ onto
$B$.  

To see this, first represent $B$ on Hilbert space, so that by (a),
$(i_A,i_G)$ becomes a covariant representation of $(A,G,\alpha)$.
It follows from the definition of $\|\cdot\|_*$ that
$i_A\times i_G$ extends to
a representation of $A\times_\alpha G$,
and then (c) implies that $i_A\times i_G$ maps $A\times_\alpha G$
onto (the represented copy of) $B$.  Now choose any
covariant representation $(\pi,U)$ of $(A,G,\alpha)$.
Then by (b), $\pi\times U(f) =
\phi_{\pi,U}(i_A\times i_G(f))$ for $f\in C_c(G,A)$, and hence for all 
$f\in A\times_\alpha G$.  But this implies $\|\pi\times U(f)\|\leq
\|i_A\times i_G(f)\|$, and hence that $\|f\|_*\leq \|i_A\times
i_G(f)\|$.  Thus $i_A\times i_G$ is injective, which proves the claim. 
\end{rem}

We next seek to define crossed products of right-Hilbert bimodules. 
There are several possible approaches; for example, we could deduce much
of the following proposition from the   construction of 
\cite[Section~5]{ComCP} for imprimitivity bimodules by first factoring the
right-Hilbert  bimodule as in 
\propref{AG-cat}(iii).  We have opted for a direct
treatment partly because it is more elementary, and partly
because we feel that the details should be available.  
This is also the approach used by Kasparov to construct crossed-product Hilbert
modules for the study of the $KK$-theory of group $C^*$-algebras in 
\cite{KasKG}[\S 6, Definition 1].  

\begin{prop}
\label{xprod}
Let $\gamma$ be an $(\alpha,\beta)$-compatible action of $G$ on a
right-Hilbert bimodule ${}_AX_B$.  There exists a right-Hilbert
$(A\times_\alpha G)$--$(B\times_\beta G)$ bimodule $X\times_\gamma G$,
which contains $C_c(G,X)$ as a dense subspace, and which satisfies
\begin{eqnarray}
\label{x1}
f\cdot h(s) & = & \int_G f(t)\cdot \gamma_t(h(t^{-1}s))\, dt,\\
\label{x2}
h\cdot g(s) & = & \int_G h(t)\cdot \beta_t(g(t^{-1}s))\, dt,
\midtext{and}\\
\label{x3}
\rip{B\times_\beta G}{h}{k}(s) & = & \int_G
\beta_{t^{-1}}(\rip{B}{h(t)}{k(ts)})\,dt
\end{eqnarray}
for $f\in C_c(G,A)$, $h,k\in C_c(G,X)$, and $g\in C_c(G,B)$.  
\end{prop}

\begin{lem}
\label{calc}
Let $\pi\times U$ be a representation of $B\times_\beta G$ on a Hilbert
space $\H$.  Then for each $h,k\in C_c(G,X)$ and $\xi,\eta\in \H$, 
\eqref{x3} satisfies
\[
(\pi\times U(\rip{B\times G}{k}{h})\xi\mid \eta)
= (\omega_{h,\xi} \mid \omega_{k,\eta}),
\]
where $\omega_{h,\xi}=\int_G h(s)\otimes U_s\xi\,ds\in X\otimes_B\H$.  
\end{lem}

\begin{proof}
We have
\begin{eqnarray*}
(\pi\times U(\rip{B\times G}{k}{h})\xi\mid \eta)
& = & \Big(\int_G \pi( \rip{B\times G}{k}{h}(s)) U_s\xi\,
ds\;\big|\; \eta\Big)\\
& = & \Big(\int_G \pi\Big( \int_G
\beta_{t^{-1}}(\rip{B}{k(t)}{h(ts)})\, 
dt\,\Big) U_s\xi\, ds\;\big|\; \eta\Big)\\
& = & \int_G\int_G (U_t^*\pi(\rip{B}{k(t)}{h(ts)})
U_{ts}\xi\mid \eta)\,ds\,dt\\
& = & \int_G\int_G (\pi(\rip{B}{k(t)}{h(s)})U_s\xi\mid
U_t\eta)\,ds\,dt\\
& = & \int_G\int_G(h(s)\otimes U_s\xi\mid k(t)\otimes
U_t\eta)\,ds\,dt\\
& = & \Big(\int_G h(s)\otimes U_s\xi\,ds\;\big|\;
\int_G k(t)\otimes U_t\eta\,dt\Big)\\
& = & (\omega_{h,\xi}\mid \omega_{k,\eta}).
\end{eqnarray*}
\end{proof}

\begin{proof}[Proof of \propref{xprod}]
We begin by showing that $C_c(G,X)$ can be
completed to give a Hilbert $B\times_\beta G$-module $X\times_\gamma G$
 satisfying
(\ref{x2}) and (\ref{x3}).  
Straightforward calculations show that (\ref{x2})
and (\ref{x3}) make $C_c(G,X)$ into a pre-inner product $B\times_\beta
G$-module, so we need only verify that the sesquilinear form of \eqref{x3}
is positive definite.  

To do so, fix, for the remainder of the proof, a faithful representation
$\pi\times U$ of $B\times_\beta G$ on a Hilbert space $\H$.  
Then for each $\xi\in\H$, \lemref{calc} gives
\[
(\pi\times U(\rip{B\times G}{h}{h})\xi\mid \xi)
= (\omega\mid \omega)\geq 0,
\]
which shows positivity.  For definiteness, suppose $h\in C_c(G,X)$
satisfies $\rip{B\times G}{h}{h} = 0$.  Then 
\lemref{calc} gives $(\pi\times U(\rip{B\times G}{h}{h})\xi\mid \xi)
= (\omega\mid \omega) =0$ for each $\xi\in\H$, so that $\omega=\int_G
h(s)\otimes U_s\xi\,ds=0$ for each $\xi$.  Thus,
\[
\Big(\int_G h(s)\otimes U_s\xi\,ds\;\big|\; x\otimes\zeta\Big) 
= \Big(\int_G \pi(\rip{B}{x}{h(s)})U_s\xi\,ds\;\big|\;\zeta\Big) 
= (\pi\times U(g_x)\xi\mid\zeta) = 0
\]
for each $\xi,\zeta\in \H$ and $x\in X$, where $g_x(s)=\rip{B}{x}{h(s)}$
defines $g_x\in C_c(G,X)$.  It follows that $g_x=0$ for all $x\in X$,
whence $h=0$ in $C_c(G,X)$.  

We next show that the $(B\times_\beta G)$-valued 
inner product on $X\times_\gamma G$ is full.
Since $C_c(G,B)$ has an approximate identity for $B\times_\beta G$,
since $\rip{B}{X}{X}$ is dense in $B$, and since $B$ acts nondegenerately
on
$B\times_\beta G$, functions of the form $f^*\rip{B}{x}{y}g$, where
$f,g\in C_c(G,B)$ and $x,y\in X$, span a dense subspace of
$B\times_\beta G$.  Now letting $(x\cdot f)(s) = x\cdot(f(s))$ define
$x\cdot f\in C_c(G,X)$, a straightforward calculation shows 
that $\rip{B\times G}{x\cdot f}{y\cdot g} = f^*\rip{B}{x}{y}g$, so that 
$\rip{B\times G}{C_c(G,X)}{C_c(G,X)}$ is dense in $B\times_\beta G$.  

We now claim that \eqref{x1} makes 
$X\times_\gamma G$ into a right-Hilbert 
$(A\times_\alpha G)$--$(B\times_\beta G)$ bimodule.  
Again, checking the algebraic conditions
of \defnref{rHb-defn} at the level of $C_c$-functions is routine; we
need to show that $\|f\cdot h\|\leq \|f\|\,\|h\|$ for $f\in C_c(G,A)$ and
$h\in C_c(G,X)$ to see that this extends 
to an action of $A\times_\alpha G$
on $X\times_\gamma G$.  

To show this, we begin by defining actions of $A$ and $G$ on $C_c(G,X)$ by
\[
(a\cdot h)(s) = a\cdot (h(s))\midtext{and} (t\cdot h)(s) =
\gamma_t(h(t^{-1}s)).
\]
Now we  use these actions to define
\[
\rho(a)(h\otimes\xi) = a\cdot h\otimes\xi
\midtext{and}
V_t(h\otimes\xi) = t\cdot h\otimes\xi;
\]
we claim that $\rho(a)$ is bounded on $C_c(G,X)\odot\H$, and hence
extends to an operator on $(X\times_\gamma G)\otimes_{B\times_\beta G}\H$.  
Again using \lemref{calc},  and writing $\omega_i$ for
$\omega_{h_i,\xi_i}$, we have
\begin{eqnarray*}
\Big\|\rho(a)\Big(\sum_{i=1}^n h_i\otimes\xi_i\Big)\Big\|^2
& = & \Big(a\cdot\sum_{i=1}^n  h_i\otimes\xi_i\mid a\cdot \sum_{j=1}^n
h_j\otimes\xi_j\Big)\\
& = & \sum_{i,j=1}^n(a\cdot h_i\otimes\xi_i\mid a\cdot h_j\otimes\xi_j)\\
& = & \sum_{i,j=1}^n (\pi\times U(\rip{B\times G}{a\cdot h_j}{a\cdot
h_i})\xi_i\mid \xi_j)\\
& = & \sum_{i,j=1}^n \Big(\int_G a\cdot h_i(s)\otimes U_s\xi_i\,ds\;\big|\; \int_G a\cdot
h_j(s)\otimes U_s\xi_j\,ds\Big)\\
& = &  \sum_{i,j=1}^n (a\cdot\omega_i\mid a\cdot\omega_j)\\
& = &  \Big(a\cdot\sum_{i=1}^n\omega_i\mid
a\cdot\sum_{j=1}^n\omega_j\Big)\\
& \leq & \|a\|^2  \sum_{i,j=1}^n (\omega_i\mid\omega_j)\\
& = & \|a\|^2  \sum_{i,j=1}^n (\pi\times U(\rip{B\times
G}{h_j}{h_i})\xi_i\mid \xi_j)\\
& = & \|a\|^2 \Big(\sum_{i=1}^n h_i\otimes\xi_i\mid  \sum_{j=1}^n
h_j\otimes\xi_j\Big)\\
& = & \|a\|^2 \Big\| \sum_{i=1}^n h_i\otimes\xi_i\Big\|^2,
\end{eqnarray*}
where the inequality holds because $A$
acts boundedly on $X\otimes_B\H$ via the induced representation
$X\dashind \pi$.

It is straightforward to check that 
each $V_t$ is unitary, and then that $(\rho,V)$ is
covariant for $(A,G,\alpha)$; $\rho$ is nondegenerate
because $A$ acts nondegenerately on $X$, so we get a nondegenerate
representation $\rho\times V$ of $A\times_\alpha G$ on
$(X\times_\gamma G)\otimes_{B\times G}\H$.  

Now for $f\in C_c(G,A)$, $h\in C_c(G,X)$,  and $\xi\in\H$, we have 
\[
\rho\times V(f)(h\otimes\xi) = \int_G \rho(f(t))V_t(h\otimes\xi)\,dt
= \Big(\int_G f(t)t\cdot h\,dt\Big)\otimes \xi
= f\cdot h\otimes\xi,
\]
so
\begin{eqnarray*}
(\pi\times U(\rip{B\times G}{f\cdot h}{f\cdot h})\xi\mid\xi)
& = & (f\cdot h\otimes\xi\mid f\cdot h\otimes\xi)\\
& = & (\rho\times V(f)(h\otimes\xi)\mid \rho\times V(f)(h\otimes\xi))\\
& \leq & \|f\|^2(h\otimes\xi\mid h\otimes\xi)\\
& = & (\pi\times U(\|f\|^2\rip{B\times G}{h}{h})\xi\mid\xi).
\end{eqnarray*}
It follows that $\|f\cdot h\|\leq\|f\|\|h\|$.  

Finally, to see that the action of $A\times_\alpha G$  on
$X\times_\gamma G$ is nondegenerate, note that for $h\in C_c(G,X)$, 
the Cauchy-Schwarz inequality gives
\begin{eqnarray*}
\|h\|^2 
& = & \|\rip{B\times G}{h}{h}\|\\
&\leq& \int_G \|\rip{B\times G}{h}{h}(s)\|\,ds\\
&\leq& \int_G\int_G \|\beta_{t^{-1}}(\rip{B}{h(t)}{h(ts)}\|\,dt\,ds\\
&\leq& \int_G\int_G \|h(t)\|\,\|h(ts)\|\,dt\,ds\\
&\leq& \Big(\int_G\|h(t)\|\,dt\Big)^2,
\end{eqnarray*}
so $\|h\|\leq\|h\|_1$.  Now standard arguments show that we can choose
$f\in C_c(G,A)$ to make $\|h-f\cdot h\|_1$ arbitrarily small: take $f$ of
the form $s\mapsto a\chi(s)$ as $a$ runs through an approximate identity
for $A$ and $\chi$ runs through an approximate identity for $C_c(G)$.  
\end{proof}

\begin{ex}\label{homos}
Suppose $(A,\alpha)$ and $(B,\beta)$ are objects in $\A(G)$, and
$\phi:A\to M(B)$ is a nondegenerate homomorphism such that
$\phi\circ\alpha_s=\beta_s\circ\phi$. Then the crossed-product morphism
$_AB_B\times G$ is the completion of $C_c(G,B)$ in the norm coming from
$B\times_\beta G$, and hence is precisely $(B\times_\beta
G)_{B\times_\beta G}$; the left action is given by a nondegenerate
homomorphism $\phi\times G:A\times_\alpha G\to M(B\times_\beta G)$.
Equation~\ref{x1} shows that for functions $f\in C_c(G,A)$ and $g\in
C_c(G,B)$, we have
\[
((\phi\times G)(f)g)(s)=\int_G \phi(f(t))\beta_t(g(t^{-1}s))\,dt.
\]
\end{ex}

\begin{prop}
\label{funct}
The maps defined by 
\[
(A,\alpha)\mapsto A\times_\alpha G\midtext{and}
[{}_{(A,\alpha)}(X,\gamma)_{(B,\beta)}]\mapsto [{}_{A\times
G}(X\times_\gamma G)_{B\times G}]
\]
give a functor from $\A(G)$ to $\C$.  
\end{prop}

\begin{proof}
We first show that the map on morphisms is well-defined.  Suppose
$\phi\colon X\to Y$ is an isomorphism of right-Hilbert $A$--$B$
bimodules which is equivariant for $(\alpha,\beta)$-compatible actions
$\gamma$ and $\eta$ of $G$.  
Then $\Phi(h)(s) = \phi(h(s))$ is easily seen to give a bijective map
$\Phi\colon C_c(G,X)\to C_c(G,Y)$ which respects the right-Hilbert
bimodule structures (Equations~\ref{x1}--\ref{x3}) and hence extends to
a right-Hilbert $(A\times_\alpha G)$--$(B\times_\beta G)$ bimodule
isomorphism of $X\times_\gamma G$ onto $Y\times_\eta G$.  
Example~\ref{homos} shows that identity morphisms go to identity
morphisms.  

It only remains to see that the assignment $[X,\gamma]\mapsto
[X\times_\gamma G]$ respects composition of morphisms; that is, if
$({}_AX_B,\gamma)$ is $(\alpha,\beta)$-compatible 
and $({}_BY_C,\eta)$ is $(\beta,\epsilon)$-compatible, we need to show
that 
\[
(X\times_\gamma G)\otimes_{B\times G}(Y\times_\eta G) \cong
(X\otimes_B Y)\times_{\gamma\otimes\eta}G
\]
as right-Hilbert $A\times_\alpha G$--$C\times_\epsilon G$ bimodules.
\footnote{An analog of this result appears as Lemma~3.10 in \cite{KasEK}, 
in the context of not-necessarily-full Hilbert modules.}

The rule
\[
\Psi(h \otimes k)(s)
= \int_G h(t) \otimes \eta_t(k(t^{-1}s)) \,dt
\]
defines a linear map from $C_c(G,X) \odot C_c(G,Y)$ to
$C_c(G,X \otimes_B Y)$ which preserves the pre-right-Hilbert bimodule
structures. In order to see that $\Psi$
extends
to an isomorphism of the completions, we need only verify that $\Psi$
has dense range for the inductive limit topology.
For this, let $x \in X$ and $f \in C_c(G,B)$, and define $h \in
C_c(G,X)$ by $h(s) = x \cdot f(s)$. Then for $k \in C_c(C,Y)$ we have
\begin{align*}
\Phi(h \otimes k)(s)
&= x \otimes \int_G f(t) \cdot \eta_t(k(t^{-1}s)) \,dt
\\&= x \otimes (f \cdot k)(s).
\end{align*}
Now, we can approximate $k$ by $f \cdot k$
in the inductive limit topology, and taking $k$ of the form $k(s) =
y g(s)$ for $y \in Y$ and $g \in C_c(G)$ we can thus approximate
the function $s \mapsto (x \otimes y) g(s)$. But such functions
have inductive-limit-dense span in $C_c(G,X \otimes_B Y)$.

\end{proof}


\section{Naturality in Green's Imprimitivity Theorem}
\label{nat}

Suppose that $\alpha$ is an action of a locally compact
group $G$ on a $C^*$-algebra $A$, and $H$ is a closed
subgroup of $G$. Takesaki showed in \cite{TakCR} how Mackey's
construction of induced representations of groups could be modified to
induce covariant representations of
$(A,H,\alpha)$ to representations of $(A,G,\alpha)$. Subsequently
Green showed how the integrated forms of Takesaki's induced
representations could be obtained using Rieffel's abstract induction by
Hilbert bimodules \cite{GreLS}. To construct his $(A\times_\alpha
H)$--$(A\times_\alpha G)$ bimodule $X_H^G(A,\alpha)$, Green completed the
$C_c(G,A)$--$C_c(H,A)$ bimodule
$C_c(G,A)$ in which
\begin{align*}
f\cdot x(t)&=\int_G f(s)\alpha_s(x(s^{-1}t))\Delta(s)^{1/2}\,ds\\
x\cdot g(t)&=\int_H x(th)\alpha_{th}(g(h^{-1}))\Delta(h)^{-1/2}\,dh\\
\langle x,y\rangle_{A\times_\alpha H}(h)&=
\int_G \alpha_s\bigl(x(s^{-1})^*y(s^{-1}h)\bigr)
\Delta(h)^{-1/2}\,ds
\end{align*}
for $x\in C_c(G,A)$, $f\in C_c(G,A)\subset A\times_\alpha G$, $g\in C_c(H,A)$.
The representation of $A\times_\alpha G$ induced from a representation
$\pi$ of $A\times_\alpha H$ is then by definition the representation
$X_H^G$-$\Ind \pi$ acting in $X_H^G\otimes_{A\times H} \H_\pi$.

Imprimitivity theorems tell us how to recognize induced 
representations. In Rieffel's modern formulation, 
we recognize a representation which has been induced from
$H$ by the presence of an auxilary representation of $C_0(G/H)$. In
Green's realization, this representation comes from the action of
$C_0(G/H)$ on $C_c(G,A)$ given by
\[
c\cdot x(t)=c(tH)x(t),
\]
and the crucial properties of this action are that it commutes with the
left action of
$A$ on $X_H^G$ and is covariant for the left action of $G$:
more formally,
\[
c\cdot(i_A(a)\cdot x)=i_A(a)\cdot (c\cdot x) \ \mbox{ and }\ 
i_G(s)\cdot (c\cdot (i_G(s)^*\cdot x))=\tau_s(c)\cdot x,
\]
where we have extended the left action of $A\times_\alpha G$ on $X_H^G$
to $M(A\times_\alpha G)$ and $\tau$ is the action of $G$ on $C_0(G/H)$
given by $\tau_s(c)(t)=c(s^{-1}t)$. 
In other words, the left actions
combine to give an action of $(A\otimes
C_0(G/H))\times_{\alpha\otimes\tau}G$ on $X_H^G(A,\alpha)$,
and it is the content of Green's Imprimitivity Theorem that
$X_H^G$ is an $((A\otimes
C_0(G/H))\times_{\alpha\otimes\tau}G)$--$(A\times_{\alpha}H)$
imprimitivity bimodule. Note that the left action
of $A\times_\alpha G$ on $X_H^G$ factors through the left action of
$(A\otimes
C_0(G/H))\times_{\alpha\otimes\tau}G$ via the canonical map of
$A\times G$ into $M((A\otimes C_0(G/H))\times_{\alpha\otimes\tau} G)$.

\begin{thm}
\label{trans-thm}
Suppose $H$ is a closed subgroup of a locally compact group $G$. Then
the assignment $(A,\alpha)\mapsto [X_H^G(A,\alpha)]$ is a natural
transformation between the functors $(A,\alpha)\mapsto A\times_\alpha G$
and 
$(A,\alpha)\mapsto A\times_\alpha H$ from $\mathcal{A}(G)$ to
$\mathcal{C}$, and a natural equivalence between the functors 
$(A,\alpha)\mapsto (A\otimes C_0(G/H))\times_{\alpha\otimes \tau} G$ and 
$(A,\alpha)\mapsto A\times_\alpha H$ from $\mathcal{A}(G)$ to
$\mathcal{C}$.
\end{thm} 

Before we begin the proof, we point out that we
have implicitly asserted in the statement of the theorem that 
the crossed products $A\times_\alpha H$ and $(A\otimes
C_0(G/H))\times_{\alpha\otimes \tau} G$ define functors on
$\mathcal{A}(G)$. The first case is quite
easy: we just need to check that restricting every action in sight
from $G$ to $H$ is a functor from $\mathcal{A}(G)$ to
$\mathcal{A}(H)$, and then compose this with the 
crossed-product functor from $\mathcal{A}(H)$ to $\mathcal{C}$. 
For the second,
we need to prove that $(A,\alpha)\mapsto (A\otimes
C_0(G/H),\alpha\otimes\tau)$ is a functor from $\mathcal{A}(G)$ to
itself, or, more generally, that for any given action $\epsilon$ 
of $G$ on
a $C^*$-algebra $C$, $(A,\alpha)\mapsto (A\otimes_{\min} C,
\alpha\otimes\epsilon)$ is a functor. This is non-trivial: a morphism
$[_AX_B,\gamma]$ goes to a morphism $[X\otimes C,\gamma\otimes\epsilon]$
based on the external tensor product $_{A\otimes C}(X\otimes
C)_{B\otimes C}$ of bimodules (see \cite[Corollary 3.38]{RW}). But all
the details are routine, so we shall omit them. 

A natural transformation $T$ between two functors $F,
G:\mathcal{A}\to\mathcal{B}$ assigns to each object $A$ of $\mathcal{A}$ a
morphism $T(A):F(A)\to G(A)$ such that, for every morphism $\phi:A\to
B$ in $\mathcal A$, the diagram
\begin{equation*}\label{G-ib-diag}
\begin{diagram}
\node{F(A)}
        \arrow{s,l}{F(\phi)}
        \arrow[2]{e,t}{T(A)}
\node[2]{G(A)}
        \arrow{s,r}{G(\phi)}\\
\node{F(B)}
        \arrow[2]{e,t}{T(B)}
\node[2]{G(B)}
\end{diagram}
\end{equation*}
commutes in $\mathcal{B}$. The transformation $T$ is a natural
equivalence if $T(A)$ is an isomorphism for all objects $A$. Green's
Theorem tells us that $[X_H^G(A,\alpha)]$ is a morphism from
$A\times_\alpha G$ to $A\times_\alpha H$ and an isomorphism from 
$(A\otimes
C_0(G/H))\times_{\alpha\otimes \tau} G$ to $A\times_\alpha H$, so we
need to show that for every morphism
$[_{(A,\alpha)}(X,\gamma)_{(B,\beta)}]$ 
in the category $\mathcal{A}(G)$, the diagrams
\begin{equation}\label{mor-diag}
\begin{diagram}
\node{A\times_\alpha G}
        \arrow{s,l}{X\times_\gamma G}
        \arrow[2]{e,t}{X_H^G(A,\alpha)}
\node[2]{A\times_\alpha H}
        \arrow{s,r}{X\times_\gamma H}\\
\node{B\times_\beta G}
        \arrow[2]{e,t}{X_H^G(B,\beta)}
\node[2]{B\times_\beta H}
\end{diagram}
\end{equation}
and
\begin{equation}\label{iso-diag}
\begin{diagram}
\node{(A\otimes C_0(G/H))\times_{\alpha\otimes\tau} G}
        \arrow{s,l}{(X\otimes C_0(G/H))\times_{\gamma\otimes\tau} G}
        \arrow[2]{e,t}{X_H^G(A,\alpha)}
\node[2]{A\times_\alpha H}
        \arrow{s,r}{X\times_\gamma H}\\
\node{(B\otimes C_0(G/H))\times_{\beta\otimes\tau} G}
        \arrow[2]{e,t}{X_H^G(B,\beta)}
\node[2]{B\times_\beta H}
\end{diagram}
\end{equation} 
commute in $\mathcal{C}$.

\begin{proof}[Proof of \thmref{trans-thm}]
We establish the commutativity of both (\ref{mor-diag}) and
(\ref{iso-diag}) by factoring 
$[X,\gamma] = [Y,\eta]\circ[\phi,\epsilon]$, 
as in \propref{AG-cat}(iii), where
$[_{(C,\epsilon)}(Y,\eta)_{(B,\beta)}]$
is an isomorphism in $\A(G)$ and 
$\phi:A\to M(C)$ is an $\alpha$--$\epsilon$ equivariant nondegenerate
homomorphism.  
We shall prove the
commutativity of (\ref{mor-diag}) by showing that we have two
commutative diagrams
\begin{equation}\label{diag1}
\begin{diagram}
\node{A\times_\alpha G}
        \arrow{s,l}{\phi\times G}
        \arrow[2]{e,t}{X_H^G(A,\alpha)}
\node[2]{A\times_\alpha H}
        \arrow{s,r}{\phi\times H}\\
\node{C\times_\epsilon G}
        \arrow[2]{e,t}{X_H^G(C,\epsilon)}
\node[2]{C\times_\epsilon H}
\end{diagram}
\end{equation}
and
\begin{equation}\label{diag2}
\begin{diagram}
\node{C\times_\epsilon G}
        \arrow{s,l}{Y\times_\eta G}
        \arrow[2]{e,t}{X_H^G(C,\epsilon)}
\node[2]{C\times_\epsilon H}
        \arrow{s,r}{Y\times_\eta H}\\
\node{B\times_\beta G}
        \arrow[2]{e,t}{X_H^G(B,\beta)}
\node[2]{B\times_\beta H.}
\end{diagram}
\end{equation}
Because taking crossed products is a functor, this will then give
\begin{align*}
[X_H^G(B)]\circ[_{A\times G}(X\times_\gamma G)_{B\times G}]&=
[X_H^G(B)]\circ[_{C\times G}(Y\times_\eta G)_{B\times G}]
\circ[\phi\times G]\\
&=[_{C\times H}(Y\times_\eta H)_{B\times H}]\circ
[\phi\times H]\circ [X_H^G(A)]\\
&=[_{A\times H}(X\times_\gamma H)_{B\times H}]\circ [X_H^G(A)],
\end{align*}
as required.

The commutativity of (\ref{diag1}) amounts to:

\begin{lem}\label{commforHMs}
Suppose $(A,\alpha)$, $(C,\epsilon)$ are objects in $\mathcal{A}(G)$ and
$\phi:A\to M(C)$ is a nondegenerate homomorphism such that $\phi\circ
\alpha=\epsilon\circ\phi$. Then 
\[
X_H^G(A,\alpha)\otimes_{A\times H}(C\times_\epsilon H)
\cong X_H^G(C,\epsilon)\cong
(C\times_\epsilon G)\otimes_{C\times G} X_H^G(C,\epsilon)
\]
as right-Hilbert $(A\times_\alpha G)$--$(C\times_\epsilon H)$ bimodules.
\end{lem} 

\begin{proof}
The second isomorphism is standard. For the first, define
\[
\Psi:C_c(G,A)\odot C_c(H,C)\to C_c(G,C)
\]
by
\[
\Psi(x\otimes
g)(t)=\int_H\phi(x(th))\epsilon_{th}(g(h^{-1}))\Delta(h)^{-1/2}\,dh.
\]
The usual change-of-variables arguments show that $\Psi$ is right
$C_c(H,C)$-linear
(though strictly speaking it is not necessary to prove this). 
The map $\phi\times H$ is  given  by
\begin{equation}\label{star}
(\phi\times H(f)g)(h)=
\int_H \phi(f(k))\epsilon_k(g(k^{-1}h))\,dk;
\end{equation}
using this identity, 
some convoluted calculations involving several
changes of variables show that
\[
\big\langle\Psi(x\otimes g),\Psi(y\otimes f)\big\rangle_{C_c(H,C)}=
\big(\phi\times H\big(\langle y,x\rangle_{C_c(H,A)}\big)g\big)^*f
\]
for $x,y\in C_c(G,A)\subset X_H^G(A)$ and $g,f\in C_c(H,C)$. Thus
$\Psi$ converts the inner product on the internal tensor product
$X_H^G(A)\otimes_{A\times H}(C\times H)$ to the usual one on
$C\times_\epsilon H$. 

To see that $\Psi$ has dense range, and hence extends to an isomorphism
of right-Hilbert bimodules, it is enough to approximate elements in
$C_c(G,C)$ of the form $y\cdot f$, where $y\in C_c(G,C)$ and $f\in
C_c(H,C)$, with elements in the range of $\Psi$.  Moreover, a partition
of unity argument shows that functions of the form $t\mapsto
\phi(x(t))\epsilon_t(c)$ for $x\in C_c(G,A)$ and $c\in C$
span an inductive-limit dense ---
and hence norm-dense --- subspace of $C_c(G,C)$,
so we can assume that $y$ has this form.
But now a routine calculation shows that if $g(h) = cf(h)$, then 
$\Psi(x\otimes g) = y\cdot f$ in $C_c(G,C)$; this implies that the
range of $\Psi$ is dense. 

Another calculation using the analogue of (\ref{star}) for the
homomorphism $\phi\times G$ shows that this isomorphism respects the
left action of $A\times_\alpha G$.
\end{proof}

To prove the commutativity of (\ref{diag2}), we use a device from
\cite[\S4]{ER-ST}. Recall that if $_AX_B$ is an imprimitivity bimodule,
then the linking algebra is the collection of $2\times 2$ matrices
\[
L(X):=\left\{
\begin{pmatrix}
a&x\\
\widetilde{z}&b\end{pmatrix}:a\in A,\ b\in B,\ x,z\in X\right\}.
\]
with multiplication and involution given by
\[
\begin{pmatrix}
a&x\\
\widetilde{z}&b\end{pmatrix}
\begin{pmatrix}
a'&x'\\
\widetilde{z'}&b'\end{pmatrix}:=
\begin{pmatrix}
aa'+{}_A\langle x,z'\rangle&a\cdot x'+x\cdot b'\\
(z\cdot a'+b\cdot z')\,\widetilde{}&\langle z,x'\rangle_B+bb'\end{pmatrix}
\]
and
\[
\begin{pmatrix}
a&x\\
\widetilde{z}&b\end{pmatrix}^*:=
\begin{pmatrix}
a^*&z\\
\widetilde{x}&b^*\end{pmatrix}.
\]
This has a unique complete $C^*$-norm, obtained, for example, by
identifying $L(X)$ with the $C^*$-algebra 
$\mathcal{K}(X\oplus B)$ of compact operators
on the Hilbert module direct sum
$(X\oplus B)_B$ (see \cite[page 50]{RW}). The matrices
\[
p=p_{L(X)}:=\begin{pmatrix}
1_{M(A)}&0\\
0&0\end{pmatrix}\ \mbox{ and }\ 
q=q_{L(X)}:=\begin{pmatrix}
0&0\\
0&1_{M(B)}\end{pmatrix}
\]
define full projections in $M(L(X))$ which allow us to 
identify $A$, $B$ and
$X$ with corners in $L(X)$; we can then use the projections $p$, $q$
to break up a Hilbert $L(X)$-module into modules over $A$ and $B$. Our
key technical result describes some relations among these submodules.
It is a mild generalization of \cite[Lemma~4.6]{ER-ST}.

\begin{prop}\label{link-prop}
Suppose $_AX_B$ and $_CY_D$ are imprimitivity bimodules, and $Z$ is a
right-Hilbert $L(X)$--$L(Y)$ bimodule. Then
\begin{itemize}
\item[(i)]
$pZp=p_{L(X)}Zp_{L(Y)}$ is a right-Hilbert $A$--$C$
bimodule;
\item[(ii)]
$qZq=q_{L(X)}Zq_{L(Y)}$ is a right-Hilbert $B$--$D$
bimodule;
\item[(iii)]
$pZq=p_{L(X)}Zq_{L(Y)}$ is a right-Hilbert $A$--$D$
bimodule;
\item[(iv)]
there is an isomorphism $\Phi:X\otimes_B qZq\to pZq$ of
right-Hilbert $A$--$D$ modules such that $\Phi(x\otimes  z)=x\cdot z$;
\item[(v)]
there is an isomorphism $\Psi:pZp\otimes_C Y\to pZq$ of
right-Hilbert $A$--$D$ modules such that $\Phi(w\otimes  y)=w\cdot y$.
\end{itemize}
\end{prop} 

\begin{rem}
Together, the last two parts of \propref{link-prop} say that the diagram
\begin{equation*}
\begin{diagram}
\node{A}
        \arrow{s,l}{X}
        \arrow[2]{e,t}{pZp}
\arrow{ese,t}{pZq}
\node[2]{C}
        \arrow{s,r}{Y}\\
\node{B}
        \arrow[2]{e,t}{qZq}
\node[2]{D}
\end{diagram}
\end{equation*}
commutes in $\C$. 
If we had such a proposition when $X$ and $Y$ are just right-Hilbert
bimodules, and the linking algebras are by definition $L(_{\K(X)}X_B)$ and
$L(_{\K(Y)}Y_D)$, we could avoid having to factor our morphisms.
However, we need to know that
$X$ and
$Y$ are imprimitivity bimodules to identify the corners $pL(X)p$ and
$pL(Y)p$ with $A$ and $C$. Since $C$ in particular occurs in the middle
of the internal tensor products in (v), it is hard to see how
this hypothesis might be avoided.
\end{rem}

\begin{lem}
\label{pqlem}
Suppose $_EW_F$ is a right-Hilbert bimodule and $P\in M(E)$, $Q\in M(F)$
are full projections. Then $PWQ$ is a right-Hilbert $PEP$--$QFQ$
bimodule.
\end{lem}

\begin{proof}
Since $PWQ$ is certainly a $PEP$--$QFQ$ submodule, and since 
\[
\langle
PWQ, PWQ\rangle_F=Q\langle
PW, PW\rangle_FQ\subset QFQ,
\]
we only have to check nondegeneracy of the
left action and fullness on the right. For nondegeneracy, we use the
fullness of $P$ to see that $PEP\cdot WQ=PEP\cdot EWQ=P(EPE)\cdot
WQ$ is dense in $PWQ$. For fullness, we use the fullness of $P$
again to see that
\begin{align*}
\langle
PWQ, PWQ\rangle_F&=Q\langle PE\cdot W, PE\cdot W\rangle_FQ\\
&=Q\langle
W, EP\cdot PE\cdot W\rangle_FQ=Q\langle
W, E\cdot W\rangle_FQ
\end{align*}
is dense in $QFQ$.
\end{proof}

\begin{proof}[Proof of Proposition~\ref{link-prop}]
For (i), apply \lemref{pqlem} with $P=p_{L(X)}\in M(L(X))$ and $Q=p_{L(Y)}\in
M(L(Y))$, and note that $A=PL(X)P$ and $C=QL(Y)Q$ because $X$ and $Y$
are imprimitivity bimodules. Parts (ii) and (iii) follow similiarly. 

For part (iv), we first note that because $(x,z)\mapsto x\cdot z$ is
bilinear, there  is a well-defined map $\Phi$ on the algebraic tensor
product
$X\odot qZq$ with the required property. We next verify  that
$\Phi$ preserves the inner product: if $x,x'\in X$ and $z,z'\in qZq$,
then the inner product $\langle x,x'\rangle_B$ is given in $L(X)$ by
$\langle x,x'\rangle_B=x^*x'$; more formally,
\[
\begin{pmatrix}
0&0\\
0&\langle x,x'\rangle_B\end{pmatrix}=
\begin{pmatrix}
0&x\\
0&0\end{pmatrix}^*
\begin{pmatrix}
0&x'\\
0&0\end{pmatrix}.
\] 
Thus
\begin{align*}
\langle x\otimes z,x'\otimes z'\rangle_D&=\langle z,\langle
x,x'\rangle_B\cdot z'\rangle_D\\
&=\langle z,(x^*x')\cdot z'\rangle_D\\
&=\langle x\cdot z,x'\cdot z'\rangle_D\\
&=\langle \Phi(x\otimes z),\Phi(x'\otimes z')\rangle_D,
\end{align*}
and $\Phi$ extends to an isometry of $X\otimes_B qZq$ into $pZq$. To see
that $\Phi$ has dense range and is therefore onto, note that $L(X)$ acts
nondegenerately on $Z$, so
\[
\range\Phi\supset pL(X)q\cdot qZq=pL(X)q\cdot q(L(X)\cdot
Z)q=p(L(X)qL(X))\cdot Zq
\]
is dense because $q$ is full.

For (v), note that the $D=qL(Y)q$-valued inner product is given by the
product in $L(Y)$, so
\begin{align*}
\langle w\otimes y,w'\otimes y'\rangle_D&=\langle y,\langle
w,w'\rangle_C\cdot y'\rangle_D\\
&=y^*\langle w,w'\rangle_{L(Y)}y'\\
&=\langle w\cdot y,w'\cdot y'\rangle_{L(Y)},
\end{align*}
which is $\langle w\cdot y,w'\cdot y'\rangle_{D}$ because $w\cdot y$ and
$w'\otimes y'$ belong to $pZq$. Thus $w\otimes y\mapsto w\cdot y$
extends to an isometry $\Psi$, as claimed, and $\Psi$ is surjective
because $p$ is full and $L(Y)$ acts nondegenerately on $Z$.
\end{proof}

\noindent
{\em Back to the proof of \thmref{trans-thm}.}
Consider the action
$\delta=\left(\begin{smallmatrix}\epsilon&\eta\\\widetilde\eta&\beta
\end{smallmatrix}\right)$ of $G$ on $L(Y)$;
by Green's theorem, $Z=X_H^G(L(Y))$ is then a right-Hilbert
$(L(Y)\times_\delta G)$--$(L(Y)\times_\delta H)$
bimodule.  To show that \eqref{diag2} commutes, we aim to apply 
\propref{link-prop} to $Z$, but this requires that we first identify
$L(Y)\times_\delta G$ with $L(Y\times_\eta G)$, and similarly for $H$.  
Now the dense subalgebra
$C_c(G,L(Y))$ consists of $2\times 2$ matrices with entries in
$C_c(G,C)$, $C_c(G,Y)$, {\em etc}., and the diagonal corners have their usual
$*$-algebraic structures as subalgebras of $C\times_\epsilon G$ and
$B\times_\beta G$. The norms are the same too: every covariant
representation of $(L(Y),\delta)$ restricts to a covariant representation
of $(B,\beta)$, and the inducing construction of \exref{rep-ex} shows that
every covariant representation of $(B,\beta)$ extends to a covariant
representation of $(L(Y),\delta)$. So we have embeddings of 
$C\times_\epsilon G$ and
$B\times_\beta G$ as corners in $L(Y)\times_\delta G$. The bimodule crossed
product $Y\times_\eta G$ also embeds: because $B\times G$ embeds
isometrically, the norm of a function 
\[
t\mapsto 
\begin{pmatrix}0&f(t)\\0&0
\end{pmatrix}\in C_c(G,L(Y))\subset L(Y)\times G
\]
is just $\|\langle f,f\rangle_{B\times G}\|$. These
embeddings combine to give an isomorphism of $L(Y\times_\eta G)$
onto $L(Y)\times_\delta G$, as desired. 
(This observation is not new: in \cite{ComCP},
Combes \emph{defines} $Y\times_\eta G$ to be the corner in
$L(Y)\times_\delta G$.)
In exactly the same way, $L(Y\times_\eta H)\cong L(Y)\times_\delta H$. 

We next have to identify the corners $pX_H^G(L(Y),\delta)p$ and
$qX_H^G(L(Y),\delta)q$ with $X_H^G(C,\epsilon)$ and $X_H^G(B,\beta)$.
But viewing
elements of $C_c(G,L(Y)\subset X_H^G(L(Y),\delta)$ 
as matrices of functions,
and similarly for $C_c(G,L(Y))\subset L(Y)\times_\delta G$ and 
$C_c(H,L(Y))\subset L(Y)\times_\delta H$, gives the desired identifications.
It now follows from Proposition~\ref{link-prop} that
\[
(Y\times_\eta G)\otimes_{B\times G}X_H^G(B,\beta)\cong
pX_H^G(L(Y),\delta)q\cong
X_H^G(C,\epsilon)\otimes_{C\times H}(Y\times_\eta H)
\]
as right-Hilbert $(C\times_\epsilon G)$--$(B\times_\beta H)$ bimodules. In
other words, (\ref{diag2}) commutes in the category $\mathcal{C}$. 

We have now established the first part of \thmref{trans-thm}, 
the natural transformation between the
functors $(A,\alpha)\mapsto A\times_\alpha G$ and $(A,\alpha)\mapsto
A\times_\alpha H$. For the second part, we have to prove that
(\ref{iso-diag}) commutes. We follow the same procedure as before: 
factor a given morphism
$[_AX_B,\gamma]$ as $[_CY_B,\eta]\circ[\phi,\epsilon]$,
and prove that two separate
diagrams commute. To see the commutativity of 
\begin{equation*}
\begin{diagram}
\node{(A\otimes C_0(G/H))\times_{\alpha\otimes\tau} G}
        \arrow{s,l}{(\phi\otimes\id)\times G}
        \arrow[2]{e,t}{X_H^G(A,\alpha)}
\node[2]{A\times_\alpha H}
        \arrow{s,r}{\phi\times H}\\
\node{(C\otimes C_0(G/H))\times_{\epsilon\otimes\tau} G}
        \arrow[2]{e,t}{X_H^G(C,\epsilon)}
\node[2]{C\times_\epsilon H,}
\end{diagram}
\end{equation*} 
we merely note that the isomorphism of Lemma~\ref{commforHMs} respects
the left action of $C_0(G/H)$, and hence is also an isomorphism of
right-Hilbert $\big((A\otimes C_0(G/H))\times_{\alpha\otimes\tau}
G\big)$--$(C\times_\epsilon H)$ bimodules. To see that
\begin{equation*}
\begin{diagram}
\node{(C\otimes C_0(G/H))\times_{\epsilon\otimes\tau} G}
        \arrow{s,l}{(Y\otimes C_0(G/H))\times_{\eta\otimes\tau} G}
        \arrow[2]{e,t}{X_H^G(C,\epsilon)}
\node[2]{C\times_\epsilon H}
        \arrow{s,r}{Y\times_\eta H}\\
\node{(B\otimes C_0(G/H))\times_{\beta\otimes\tau} G}
        \arrow[2]{e,t}{X_H^G(B,\beta)}
\node[2]{B\times_\beta H}
\end{diagram}
\end{equation*}
commutes, we again apply  Proposition~\ref{link-prop} to 
$Z=X_H^G(L(Y),\delta)$,
this time viewed as a $L((Y\otimes
C_0(G/H))\times_{\eta\otimes\tau}G)$--$L(Y\times_\eta H)$
imprimitivity bimodule.  To do so, we have to identify 
$L((Y\otimes C_0(G/H))\times_{\eta\otimes\tau}G)$ with 
$(L(Y)\otimes C_0(G/H))\times_{\delta\otimes\tau}G$,
but there is no essential difference in the argument.

This concludes the proof of our main theorem.
\end{proof}

\begin{rem}
As an example of how the present set of ideas can simplify
things, notice that the commutativity of diagram (\ref{diag2}) implies
that for each covariant representation $\rho\times V$ of
$B\times_\beta H$, the representation $(Y\times_\eta
G)\dashind(\Ind_H^G (\rho\times V))$ is equivalent to
$\Ind_H^G((Y\times_\eta H)\dashind(\rho\times V))$.  This is the
content of \cite[Theorem~3]{EchME}, or at least the content of its
proof (the twists in that theorem go along for free by Remark~(2)
on page~174 of \cite{EchME} and \cite[Corollary~5]{GreLS}).  The
point is that the arguments we give here are more elegant as well
as more powerful. 
\end{rem}




\end{document}